\documentclass[12pt,a4paper]{article}
\usepackage{amsthm,amsmath,amssymb,amsfonts,color, enumerate}
\usepackage[authoryear]{natbib}
 \usepackage{tikz}
\usepackage[colorlinks, allcolors=blue]{hyperref}
\usepackage{graphicx}
\usepackage{array}
\usepackage{float}
\usepackage{titling}
\usepackage{multirow}

\numberwithin{equation}{section}

\newtheoremstyle{dotless}{}{}{\itshape}{}{\bfseries}{}{ }{}
\theoremstyle{dotless}
\newtheorem{theorem}{Theorem}[section]

\theoremstyle{definition}

\newtheoremstyle{definition}{}{}{}{}{\bfseries}{}{ }{}
\theoremstyle{definition}

\newtheorem{example}{Example}[section]

\newcommand{\R}{\mathbb{R}}

\newcommand{\bea}{\begin{eqnarray*}}
\newcommand{\eea}{\end{eqnarray*}}
\newcommand{\be}{\begin{eqnarray}}
\newcommand{\ee}{\end{eqnarray}}

\newcommand{\ww}{\omega}

\usepackage[authoryear]{natbib}
\begin{document}

\title{Optimal designs for estimating individual coefficients in polynomial regression with no intercept}
\author{\small Holger Dette \\
\small Ruhr-Universit\"at Bochum \\
\small Fakult\"at f\"ur Mathematik \\
\small 44780 Bochum, Germany \\
\small e-mail: holger.dette@rub.de\\
\and
\small Viatcheslav B. Melas \\
\small St. Petersburg State University \\
\small Department of Mathematics \\
\small St. Petersburg ,  Russia \\
{\small email: vbmelas@yandex.ru}\\
\and
\small Petr Shpilev\\
\small St. Petersburg State University \\
\small Department of Mathematics \\
\small St. Petersburg , Russia \\
{\small email: pitshp@hotmail.com}\\
}

\maketitle

\begin{abstract}
In a seminal paper \cite{studden1968} characterized  $c$-optimal designs in regression models, where the regression functions form a Chebyshev system.
He used these results  to  determine the optimal design for estimating the individual coefficients in a polynomial regression model on the interval $[-1,1]$
explicitly. In this note
we identify  the  optimal design for estimating the individual coefficients in a polynomial regression model with no intercept (here the regression functions do not
form a Chebyshev system).
 \end{abstract}



\section{Introduction} \label{sec1}
\def\theequation{1.\arabic{equation}}
\setcounter{equation}{0}

Consider the common polynomial regression model of degree $n$ with no intercept
\be
 \label{1.1} Y_i = (x_{i},x_{i}^{2}, \dots, x_{i}^n)^\top \theta  + \varepsilon_i , \qquad i=1,\ldots,N,
\ee
where $\varepsilon_1,\dots,\varepsilon_N$ denote independent  random variables with $\mathbb{E}[\varepsilon_i]
=0;$ $ {\rm Var}(\varepsilon_i)=\sigma^2 >0$ $(i=1,\dots,N)$,
$\theta = (\theta_{1}, \ldots , \theta_{n})^{\top} \in \mathbb{R}^n$ is a vector of unknown parameters and the explanatory variables
$x_{1}, \ldots , x_{N}$ vary in the interval $[-1,1]$.
An (approximate) optimal  design
minimizes  an appropriate functional of the (asymptotic) covariance matrix of  the statistic  $\sqrt{N} \hat \theta $,
where the  $  \hat \theta $  denotes  the least squares estimate of the parameter $\theta$
in the regression model \eqref{1.1} [see \cite{silvey1980} or \cite{pukelsheim2006}]. Numerous authors have worked  on the
problem of determing optimal designs in this model, where the main focus is on the $D$- and $E$-optimality criterion
corresponding to the minimization of the determinant  and maximum eigenvalue of the (asymptotic) covariance matrix of the least squares estimate
[see  \cite{huangchangwong1995,changheil1996,ortizrod1998,chang1999,fang2002} or \cite{lilauzhang2005}].
While these  problems have  been nowadays  well understood there exist basically no solutions
of the  optimal design problem for  other type of optimality criteria.

In the present note we add to this literature and determine explicitly the approximate (in the sense of \cite{kiefer1974}) optimal design for estimating the individual
coefficients in  a polynomial regression model with no intercept on the interval $[-1,1$]. The corresponding optimality criteria are special cases of the well known
$c$-optimality   criterion which seeks for a design minimizing the variance of the best linear unbiased estimate of the linear combination $c^\top\theta$
in model \eqref{1.1}, where $c \in \R^n$ is a given vector.
In a seminal paper \cite{studden1968}  characterizes   $c$-optimal designs  in regression models
with regression  functions forming  a Chebyshev system. As an application  he found the optimal designs
for estimating the individual coefficients in a regression with intercept, that is  $ Y_i = \sum_{\ell=0}^n \theta_{\ell} x_{i}^\ell + \varepsilon_i $.
It is also indicated in \cite{studden1968}
that in general the solution of the $c$-optimal design problem is an extremely difficult one, in particular
if  the regressions functions do not form a Chebyshev system, such as in    model \eqref{1.1}, if the explanatory variable varies int he interval $[-1,1]$.

 In Section \ref{sec2} we introduce  the basic optimal design problem and
review a geometric characterization of $c$-optimal designs. The main result can be found in Section \ref{sec3}
where the optimal designs for estimating the individual coefficients in polynomial regression model
with no intercept are determined explicitly and the theory is illustrated by several examples.

\section{$c$-optimal designs} \label{sec2}
\def\theequation{2.\arabic{equation}}
\setcounter{equation}{0}

Following \cite{kiefer1974} we call a  probability
measure
\be
 \label{1.3}
\xi=
\begin{pmatrix}
x_1 & x_2 & \cdots & x_{m}\\
\ww_1 & \ww_2 & \cdots & \ww_{m}\\
\end{pmatrix}\;
\ee
with finite support $x_1,\ldots,x_m \in [-1, 1]$ and  corresponding weights $\ww_1,\ldots,\ww_m$
 an approximate design on the  interval $[-1, 1]$.  We define
\be
 \label{1.2}
f(x) = (x, \dots, x^n)^\top
\ee
as the vector of regression functions  in the polynomial regression  model \eqref{1.1}, and   by
 $$
M(\xi) =\int_ {-1}^1f(x)f^\top (x)\xi(dx)
 $$
 the information matrix  of the design $\xi$. The interpretation of $\xi$  and  $M(\xi)$ is as follows.
 If an experimenter takes $n_{1}, \ldots  , n_{m} $ observations at the experimental  conditions $x_{1}, \ldots  , x_{m} $,
 respectively, $N= \sum_{i=1}^m n_{i}$ denotes the total sample size and $n_{i}/N$ converge to $\ww_{i}$ ($i=1, \ldots , m$),
 then the asymptotic covariance matrix of the scaled least squares estimate $\sqrt{N} \hat \theta$
  in the regression model \eqref{1.1} is given
 by $\sigma^{2}M^{-1}(\xi)$, where $\sigma^{2}$ is the variance of the  errors. An approximate  optimal design minimizes
 a functional of the matrix $M^{-1}(\xi)$ (or more generally of a generalized inverse $M^{-}(\xi)$),  which is called optimality criterion 
 in the literature  [see \cite{silvey1980} or \cite{pukelsheim2006}].

 In this paper we investigate a special case of the $c$-optimality criterion, which is defined by
  \be \Phi_{c} (\xi)=
  \begin{cases}
    c^{\top}M^{-}(\xi)c & \text{ if there exists  a vector} v \in \mathbb{R}^{n}    \text{  such that  } c = M(\xi)v;, \\
    \infty, & \text{otherwise}
  \end{cases}
  \label{1.4}
    \ee
  for a given vector $c \in \R^{n}$.  In the first case the
 design $\xi $ is called {\it admissible  for estimating the linear combination $c^\top \theta$ }  in the regression model \eqref{1.1} and the
  value of the quadratic form does not depend on the choice of the generalized inverse
 [see \cite{pukelsheim2006}].
 The criterion \eqref{1.4} corresponds to the minimization of the asymptotic variance of
  the best linear unbiased estimate for the linear combination $c^\top \theta$. In particular for the $p$th unit
vector    $e_p=(0,\ldots,0,1,0,\ldots,0)^\top \in \R^{n} $ we obtain
$ e_p ^\top \theta = \theta_{p} $  and the $e_{p}$-optimal design minimizes the asymptotic  variance
of the best linear unbiased estimate for the coefficient  $\theta_{p}$ corresponding to the monomial $x^p$ in the polynomial regression
model with no intercept  ($p=1, \ldots , n$).  Throughout this paper we denote the optimal design with respect to the criterion
$ \Phi_{e_{p}} $, which is obtained  from \eqref{1.4} for $c=e_{p}$  as $e_{p}$-optimal design or optimal design for estimating the
coefficient $\theta_{p}$ in the polynomial regression model with no intercept.

We conclude this section  with a geometric characterization  of $c$-optimal designs  called Elfving's theorem
   [see \cite{elfving1952}], which will be used in Section \ref{sec3}. A proof can be found in  \cite{DETTE2004201}.
\begin{theorem} \label{Elfving}
An admissible  design $\xi^{*}$  for estimating the linear combination $c^\top \theta$  with support points $x_1,  x_2,   \ldots , x_{m} $ and weights
$\omega_1 , \omega_2,  \ldots , \omega_{m} $  is $c$-optimal if  and only if there exists a
vector   $u \in \mathbb{R}^{d}$  and a constant $h$ such that the following conditions are satisfied:
\begin{itemize}
	\item[(1)]
	$|u^\top f(x)|\leq1$ for all $ x\in\mathcal{X}$;
	\item[(2)]
	 $|u^\top f(x_i)   |=  1 $ for all $ i=1,2,\ldots,m$ ;
	\item[(3)]
	$c = h\sum_{i=1}^{m} f(x_i)\omega_i u^\top f(x_i)$.
\end{itemize}
Moreover, in this case we have $c^\top M^{-}(\xi^*)c=h^2.$
\end{theorem}

\section{Optimal designs for  estimating individual  coefficients in models with no intercept}
\label{sec3}
\def\theequation{3.\arabic{equation}}
\setcounter{equation}{0}

For the polynomial regression model with no intercept the function $u^{\top} f$ in Theorem \ref{Elfving} is of the form
$u^\top f(x) = \sum_{\ell =1}^{n} b_{\ell} x^{\ell}$. This function will be called extremal polynomial throughout this paper.
From Theorem \ref{Elfving}  it follows that  the support points of the $e_{p}$-optimal design are the extremal
points of a - in some sense - optimal polynomial. In fact it is possible
to identify these optimal polynomials explicitly. For this purpose let
$$
T_{s}(x) =\cos (s\arccos (x) )
$$
denote  the $s$th Chebyshev polynomial of the first kind [see  \cite{szego1975}] and consider
the polynomials
\be
\label{2.1}
  T_{2k-1}(x) ~, ~~T_{2k+1}(x)
\ee
and the polynomial
\be
\label{2.2}
  E_{2k}(x) = T_k\Big ((x^2(1+\cos\frac{\pi}{2k})-\cos\frac{\pi}{2k})\Big ).
  \ee
 It is easy to see that $T_{2k-1}$ and $T_{2k+1}$ have exactly  $2k$ and $2k+2$ extremal points, which are denoted by
$s_1<s_2<\ldots<s_{2k}$ and  $x_1<x_2<\ldots<x_{2k+2}$, respectively. Note that  these points are given explicitly by
   \begin{equation}\label{chebypoints}
  s_{i} = \cos \big ( \tfrac{(2k-i)\pi}{n} \big )   ~~(i=1,2,\dots, 2k),   ~~~~~
    x_{i} = \cos \big ( \tfrac{(2k+2-i)\pi}{2k+1} \big )  ~~ (i=1,2,\dots, 2k+2).
   \end{equation}
Similarly, the polynomial $E_{2k}$ in \eqref{2.2}
has $2k$ extremal points $t_{1}, \ldots , t_{2k}$, which are given by
\begin{equation}
\label{2.1}
t_{i} = - \sqrt{\dfrac{\cos\frac{(i-1) \pi}{k}+\cos\frac{\pi}{2k}}{1+\cos\frac{\pi}{2k}}}~, ~t_{2k+1-i} = \sqrt{\dfrac{\cos\frac{(i-1) \pi}{k}+\cos\frac{\pi}{2k}}{1+\cos\frac{\pi}{2k}}}~ ,\ i=1,\ldots,k
\end{equation}
Finally for a given set of support points of a design, say  $t_1^*, \ldots , t_{m}^*$,  we define for $i=1,\ldots,m$
\begin{equation}\label{3.2a}
 \bar L_i(x)=\dfrac{x\prod_{j\neq i}(x-t_j^{*})}{t_i^{*}\prod_{j\neq i}(t_i^{*}-t_j^{*})}
\end{equation}
as  the $i$th  Lagrange basis  interpolation polynomial  without intercept  corresponding to  the nodes
 $t_{1}^{*}, \ldots ,  t_{m}^{*}$ (note that the degree of $ \bar L_i(x)$ is $m$). The main result of this paper is the following.

\begin{theorem} \label{Theorem2.1}
Consider the polynomial regression model of degree $n \geq 1$ with no intercept.
\begin{itemize}
\item[(a)] If $n=2k+1$ or $n=2k$ for  some $k\geq 1$ and $p$ is even, then there exists an
$e_p$-optimal design  supported at the extremal points $t_1,\ldots,t_{2k}$ of the polynomial $E_{2k}(x)$ defined in
\eqref{2.1}.
\item[(b)] If $n=2k$ and p is odd,  then there exists an $e_p$-optimal design supported at 
 the extremal points   $s_1,\ldots,s_{2k}$ of the polynomial $T_{2k-1}(x)$ defined in \eqref{chebypoints}.
\item[(c)] If $n=2k+1$ and $p=1$ then there exist  {exactly}  two $e_p$-optimal designs with $2k +1$ support points: one design with  support
$x_2,\ldots,x_{2k+2}$ and the other  design with support points $x_1,\ldots,x_{2k+1}$. \\
 If $n=2k+1$ and $p$ is odd, $p>1$  then there exist  {exactly}   two $e_p$-optimal designs with $2k+1$ support points. One
 design  with  support points  $x_1,\ldots, x_k, x_{k+2}\ldots, x_{2k+2}$ and the other design
  with  support points   $x_1,\ldots, x_{k+1}, x_{k+3}\ldots, x_{2k+2}$.
\end{itemize}
The weights $\ww_{1}, \ldots , \ww_{m}$ at the support points $t_{1}^{*} , \ldots , t_{m}^{*}$ of the $e_p$-optimal design are given by the formula
\begin{equation}  \label{weightder}
\omega_i= \dfrac{|a_{i,p}|}{\sum_{j=1}^{m}|a_{j,p}|}~,~ i=1,\ldots,m,
\end{equation}
where  $m=2k$ in  cases (a) and (b),  $m=2k+1$ in case (c)  and  $a_{p,i}$ is  the coefficient
of the monomial  $x^p$ in the polynomial $\bar L_i$ defined  in (\ref{3.2a}) $(i= 1, \ldots ,m)$.
\end{theorem}

{\bf Proof.} We first consider assertion $(a)$ and use Theorem \ref{Elfving} with the polynomial $u^{\top} f(x) = E_{2k} (x) $ defined in \eqref{2.2}.
The properties (1) and (2) are obviously fulfilled and it remains to show that  condition (3) holds for  some nonnegative weights $\ww_i$,
$i=1,2,\ldots,2k$. This condition reads as follows
\be \label{h1}
\delta_{qp} = h \sum^{2k}_{i=1} t^q_i \omega_i E_{2k}(t_i)~,~~q=1, \ldots , 2k+1,
\ee
where $\delta_{qp}$ denotes Kronecker's symbol.
We show that a solution  is in fact possible under the symmetry assumption $\ww_{2k-i+1}=\ww_i$, $i=1,2,\ldots,k.$ Observing that
\be \label{sym1}
&&E_{2k}(t_i)=E_{2k}(t_{2k-i+1}), \\
&& t_i^{2q+1}=-\left(t_{2k-i+1}\right)^{2q+1},\ q=0,1,\ldots,k
\label{sym2}
\ee
we  see that the condition \eqref{h1} is obviously satisfied for odd exponents (note that $p$ is even)
Consequently,  it remains to show  that there exist nonnegative weights $\ww_1,\ldots,w_{2k}$ such that
\bea
h\sum_{i=1}^{2k} t_i^{2q} \ww_i E_{2k}(t_i)=\delta_{2q,p},
\eea
which reduces using the symmetries in \eqref{sym1} and \eqref{sym2} to 
\be \label{h10}
h\sum_{i=1}^{k} t_i^{2q}  \ww_i E_{2k}(t_i)=\frac{1}{2}\delta_{2q,2p} ,~~q=1, \ldots , k
\ee
for some constant $h$ .

For this purpose we introduce the notation $\tilde{\beta}=(\beta_1,\ldots,\beta_{k})^\top$, where $\beta_i=h\ww_i E_{2k}(t_i) $, and  
$\tilde{e}_{p/2}=$ $(0,\ldots,0,1/2,0,\ldots,0)^\top\in \mathbb{R}^k$, where 1/2 is in the $p/2$ position (recall  that $p$ is even)
and  rewrite the equations  in \eqref{h10} as follows
\bea
F\tilde{\beta}=\tilde{e}_{p/2},\ 
\eea
where the matrix $F$ is defined by  $F=\left(t_i^{2q} \right)_{q,i=1}^{k}$.
Because   the functions $t^2,t^4,\ldots,t^{2k}$ generate a Chebyshev system on the interval $(-1,0),$  the matrix $F$ is non-singular and the elements of $F^{-1}$ are alternating in sign. Consequently, the components of the vector
\bea
\tilde{\beta}=F^{-1}\widetilde{e}_{p/2}
\eea
are also alternating in sign and the corresponding weights $\omega_i =  \beta_i/(hE_{2k}(t_i))$ are positive, which completes the proof of assertion (a).

\medskip

Next we consider assertion (b) , where $n=2k$ and $p$ is odd.  A direct calculation shows that properties (1) and (2)
are fulfilled for the polynomial $u^\top f(x)=T_{2k-1}(x).$ Again we have to prove the existence of nonnegative weights $\ww_i,$ $i=1,\ldots,2k$ satisfying  
part  (3) of Theorem \ref{Elfving}. We consider first the equations corresponding to even exponents and note that for arbitrary $\ww_j,$ $i=1,\ldots,2k,$ satisfying $\ww_{2k-i+1}=\ww_{i},$ $i=1,\ldots,k$ we have
\bea
\sum_{i=1}^{2k} s_i^{2q} \ww_i T_{2k-1}(s_i)=0,\ \ q=1,\ldots,k,
\eea
where we used the symmetry properties
\bea
   T_{2k}(s_i)=-T_{2k-1}(s_{2k-i+1}),\ s_i^{2q}=(s_{2k-i+1})^{2q}~,\ q=0,\ldots,k.
\eea
Therefore it  remains to consider the equations corresponding to odd exponents, i.e. there exist nonnegative weights $\ww_i,\ldots,\ww_{2k}$ such that $\ww_{i}=\ww_{2k-i+1},$ $i=1,\ldots,k$ and
\bea
h\sum_{i=1}^{2k}s_i^{2q-1}  \ww_i T_{2k-1}(s_i)  =\delta_{2q-1,p} ,\ \ q=1,\ldots,k,
\eea
 which reduce (observing the symmetry properties) to
\bea
h\sum_{i=1}^{k} s_i^{2q-1}  \ww_i T_{2k-1}(s_i)  =\frac{1}{2}\delta_{2q-1,p}
\eea
for some nonnegative $\ww_i$, $i=1,\ldots,k,$. With the notation $\tilde{\beta}=(\tilde\beta_1,\ldots,\tilde\beta_{k}),$ where $h\tilde\beta_i=\ww_i T_{2k-1}(s_i),$
and  
$\tilde{e}_{(p-1)/2}=(0,\ldots,0,1/2,0,\ldots,0)^\top\in \mathbb{R}^k$, where the non-vanishing entry  1/2 is in the $(p-1)/2$ position,  
we   rewrite these equations in matrix form
\bea
F\tilde{\beta}=\tilde{e}_{(p-1)/2},
\eea
where $F=\left( s_i^{2q-1}\right)_{q,i=1}^{k}$.
Note that   the functions $t,t^3,\dots,t^{2k-1}$ generate a Chebyshev system on the interval $(-1,0)$. Consequently, the matrix $F$ is non-singular and the elements of $F^{-1}$ are alternating in sign. This implies that the components of the vector
\bea
\tilde{\beta}=F^{-1}\tilde{e}_{(p-1)/2}
\eea
are also alternating in sign and  the corresponding weights $\omega_i =  \beta_i /(h T_{2k-1}(s_{2i-1}))$ are positive.

\medskip

In order to prove part (c)  we use  the polynomial $u^\top f(x) = T_{2k+1}(x)$ as  an extremal polynomial in Theorem \ref{Elfving} as it satisfies conditions (1) and (2) of this theorem.
Consequently, the points $x_1,\ldots,x_{2k+2}$ in \eqref{chebypoints}   are potential support points of the $e_{p}$-optimal design.  
We now choose $2k+1$ points   $t_1^*,  t_2^*,   \ldots , t_{2k+1}^* $  from the extremal  points  as described in part (c) 
of Theorem \ref{Theorem2.1}.

By Theorem  \ref{Elfving}  a design  with weights
$\omega_1 , \omega_2,  \ldots , \omega_{2k+1}$  at  the points   $t_1^*,  t_2^*,   \ldots , t_{2k+1}^* $ is  $e_p$-optimal
if  
\begin{equation} \label{h0} 
e_p=hF\beta,
\end{equation} 
for some  constant  $h$, where $\beta$  is a $(2k+1)$-dimensional  vector with components $ \beta_i=u^\top f(t_i^*)\ww_i= T_{2k+1} (t_i^*)\ww_i $ ($i=1, \ldots , 2k+1$)
 and  $F=(f(t_1^*), \ldots , f(t_{2k+1}^*))$.
Observing   the identity  $ F^{-1}F =I_{2k+1}$ (here $I_{2k+1}$ is the identity   matrix) it follows 
$$ 
e_i^\top F^{-1}f(t_j^{*})=\delta_{ij}   ~~~(i,j=1, \ldots , 2k+1) .
$$
 As these equations characterize the $i$th basis Lagrange interpolation polynomial
with knots
$t_{1}^{*}, \ldots , t_{2k+1}^{*}$ we have for any point $z \in \R$
  $$e_i^\top F^{-1}f(z)= \bar L_i(z) =   a_{i}^{\top} f(z) ~ ,~~i=1,\ldots,2k+1,$$ 
  where
\begin{equation} \label{basic} 
a_{i}  = (F^{-1})^\top  e_{i} =(a_{i,1}, \ldots , a_{i,2k+1})^{\top}
\end{equation}  
is the vector of coefficients of the $i$th basis Lagrange interpolation polynomial $(i=1, \ldots , 2k+1)$. 
Therefore we obtain for the solution of \eqref{h0}  
$$
h \beta= F^{-1} e_{p}= (a_{1,p},\ldots,a_{2k+1,p})^T 
$$
or equivalently (since $ {\beta}_i = \omega_i T_{2k+1}(t_i^{*})$)
\begin{equation}\label{cond4}
h  \beta_{i} = h \omega_{i}  T_{2k+1}(t_i^{*})  = {1\over p!} {d^{p} \over d^{p} z} 
\bar L_i(z) \Big |_{z=0 } = a_{i,p}~,~~i=1, \ldots , 2k+1.
\end{equation}

Therefore   the representation \eqref{weightder} follows if 
 $T_{2k+1} (t_1^*)a_{1,p} , \ldots ,T_{2k+1} (t_{2k+1}^*)a_{2k+1,p} $ have the same sign. 
 In this case part  (3) of Theorem \ref{Elfving} is also satisfied (as we can solve  \eqref{h0}
 with positive weights) and the part (c) of Theorem \ref{Theorem2.1} proved.
For a proof of this property  we  now consider the different cases in Theorem \ref{Theorem2.1} separately. 
\\
First consider the case $p=1$ and   let  $t_1^*,\dots,t_{2k+1}^*$ be either $x_1,\ldots,x_{2k+1}$ or $x_2,\ldots,x_{2k+2}$. 
Note that in this case  either the smallest point $-1$ or the largest point  $1$ has been  deleted from the whole set of the extremal points of 
the Chebyshev polynomial $T_{2k+1}(x)$.      A direct calculation by Vieta' formulas gives for the 
$i$th  coefficient  of  the polynomial  \eqref{3.2a} 
\bea
a_{i,1}=\frac{\prod_{j=1}^{2k+1} t_j^*}{(t_i^*)^2\prod_{j\neq i}(t_i^*-t_j^*)}~,~~ i=1,\dots,2k+1,
\eea
(note that the polynomial $\bar L_i(z) =a_{i}^{T} f(z) $ in  \eqref{3.2a}  has the roots $t_1^*,\dots,t_{2k+1}^*$ and $0$).
As  the sign of the denominator is alternating with $i$ and   the sign of  $T_{2k+1}(t_i^*)$  
is also alternating with $i$  it follows that all products   $T_{2k+1} (t_i^*)a_{i,1}$ have the same sign, $i=1,2,\dots,2k+1$
(note that the numerator does not depend on $i$). 
\\
In the case  where   $p=2l+1>1$ is odd  the argument is very similar.  Here 
let $t_1^*,\dots,t_{2k+1}^*$ be either $x_1,x_2, \ldots,x_k, x_{k+2},\ldots,x_{2k+2}$ or $x_1,x_2\ldots,x_{k+1}, x_{k+3},\ldots,x_{2k+2}$.
This means  that in this case  one of the  two points  with minimal distance to $0$  has been deleted from the  set of the extremal points of 
$T_{2k+1}(x)$.  By  the  Vieta' formulas we obtain for the $i$th coefficient of the polynomial 
 $\bar {L}_{2l+1} (z) $ in   \eqref{3.2a} the representation 
\bea
a_{i,2l+1}=-\frac{\sum_{\substack{ 1\leq j_1<j_2<\ldots<j_{2l}\leq 2k+1 \\ j_{1} , \ldots  , j_{2l} \neq i}} \prod_{s=1}^{2l} t_{j_s}^*}{t_i^*\prod_{j\neq i}(t_i^*-t_j^*)}~,~~i=1,\dots,2k+1
\eea
 (note that one of the roots is equal to $0$)
and the symmetry of the roots yields
\bea
a_{i,2l+1}=-\frac{\sum_{\substack{ 1\leq j_1<j_2<\ldots,j_l\leq k+1 \\ j_{1} , \ldots  , j_{l} 
 \not  \in \{ i,  2k+2-i \}} }\prod_{s=1}^l (t_{j_s}^*)^2}{t_i^*\prod_{j\neq i}(t_i^*-t_j^*)}~,~i=1,\dots,2k+1.
\eea
Now it can be easily checked that 
 $T_{2k+1}(t_1^*)a_{1,2l+1}, \ldots T_{2k+1}(t_{2k+1}^*)a_{2k+1,2p+1}$ have the same sign.
These arguments  complete the proof of part  (c) of Theorem \ref{Theorem2.1}.

\medskip
Finally, it remains to show   the representation  \eqref{weightder} for the weights   in the  case (a) and (b). 
We omitt the details here as  this can be done in a similar way as in the proof of  part  (c) of Theorem \ref{Theorem2.1}.
\hfill $\Box$

\bigskip
\begin{example} \label{ex1}
{\rm We determine the optimal designs for estimating  the  individual coefficients in a cubic regression
 with no intercept. For this purpose let $P (x)$ be an extremal polynomial from Elfving's theorem.
 \begin{itemize}
\item[(a)] If  $p = 1$  we can use part (c) of  Theorem \ref{Theorem2.1}.
 The extremal polynomial
is given by $P(x)=x^3-\frac{3}{4}x$ with extremal points  $-1 $, $-\frac{1}{2}$, $  \frac{1}{2}$
and $ 1$. There exist two $3$-point $e_{1}$-optimal designs. One with masses $ \frac{1}{9}$, $ \frac{2}{3}$ and $ \frac{2}{9}$
at the points $-1 $, $-\frac{1}{2}$, and  $  \frac{1}{2}$ and the other one
with masses $ \frac{2}{9}$, $ \frac{2}{3}$ and $ \frac{1}{9}$ at the points  $-\frac{1}{2}$, $  \frac{1}{2}$  and $ 1$.
\item[(b)] If  $p = 2$  we can use part (a) of  Theorem \ref{Theorem2.1}.  Consequently, there exists a  unique $e_{2}$-optimal
design supported at $2$ points, that is
\bea
\begin{pmatrix}
-1 & 1\\
\frac{1}{2} & \frac{1}{2}\\
\end{pmatrix}.
\eea
In this case the corresponding extremal polynomial is not unique and given by
$ P(x)=x^2-qx+qx^3$, where $q\in[-1,1]$.
\item[(c)] If  $p = 3$  we can again use   part (c) of  Theorem \ref{Theorem2.1}.
The extremal polynomial
is given by $ P(x)=x^3-\frac{3}{4}x$
with extremal points  $-1 $, $-\frac{1}{2}$, $  \frac{1}{2}$
and $ 1$. There exist two $3$-point $e_{3}$-optimal  designs. One with masses $ \frac{1}{12}$, $ \frac{2}{3}$ and $ \frac{1}{4}$
at the points $-1 $, $ \frac{1}{2}$, and  $ 1$ and the other one
with masses $ \frac{1}{4}$, $ \frac{2}{3}$ and $ \frac{1}{12}$ at the points  $-1$, $-\frac{1}{2}$  and $ 1$.
\end{itemize}
}
\end{example}

\begin{example} \label{ex2}
{\rm We determine the optimal designs for estimating  the  individual coefficients in  a  polynomial regression model
of degree four  with no intercept. 
Note that in this case  Theorem  \ref{Theorem2.1}(a)  for $p=2,4$  and Theorem  \ref{Theorem2.1}(b)  for $p=1,3$ are applicable. Consequently the $e_{p}$-optimal designs are always
unique
\begin{itemize}
\item[(a1)]   If $p = 2$, the extremal polynomial is given by $P(x)=x^4-2(\sqrt{2}-1)x^2$
and the  unique $4$-point optimal design for estimating the coefficient of $x^{2}$
is given by
\bea
\begin{pmatrix}
-1 & -\sqrt{\sqrt{2}-1} & \sqrt{\sqrt{2}-1} & 1\\
\frac{\sqrt{2}}{8\sqrt{2}+8} & \frac{3\sqrt{2}+4}{8\sqrt{2}+8} & \frac{3\sqrt{2}+4}{8\sqrt{2}+8} & \frac{\sqrt{2}}{8\sqrt{2}+8}\\
\end{pmatrix} .
\eea
\item[(a2)]   If $p = 4$, the extremal polynomial is given by $P(x)=x^4-2(\sqrt{2}-1)x^2$
and the  unique $4$-point optimal design for estimating the coefficient of $x^{4}$
is given by
\bea
\begin{pmatrix}
-1 & -\sqrt{\sqrt{2}-1} & \sqrt{\sqrt{2}-1} & 1\\
\frac{\sqrt{2}}{4\sqrt{2}+4} & \frac{\sqrt{2}+2}{4\sqrt{2}+4} & \frac{\sqrt{2}+2}{4\sqrt{2}+4} & \frac{\sqrt{2}}{4\sqrt{2}+4}\\
\end{pmatrix} .
\eea
\item[(b1)]   If $p = 1$ , the extremal polynomial is given by $P(x)=x^3-\frac{3}{4}x$
and the  unique $4$-point optimal design for estimating the coefficient of $x^{1}$
is given by
\bea
\begin{pmatrix}
-1 & -\frac{1}{2} & \frac{1}{2} & 1\\
\frac{1}{18} & \frac{4}{9} & \frac{4}{9} & \frac{1}{18}\\
\end{pmatrix} .
\eea
\item[(b2)]   If $p = 3$ , the extremal polynomial is given by $P(x)=x^3-\frac{3}{4}x$
and the  unique $4$-point optimal design for estimating the coefficient of $x^{3}$
is given by
\bea
\begin{pmatrix}
-1 & -\frac{1}{2} & \frac{1}{2} & 1\\
\frac{1}{6} & \frac{1}{3} & \frac{1}{3} & \frac{1}{6}\\
\end{pmatrix} .
\eea
\end{itemize}
}
\end{example}

Note that this design is also optimal for estimating the coefficient of $x^3$ and in  a cubic  regression with intercept [see \cite{dette1990}].

\medskip\medskip\medskip\medskip

\noindent
{\bf Acknowledgements}
This work has been supported in part by the
Collaborative Research Center ``Statistical modeling of nonlinear
dynamic processes'' (SFB 823, Teilprojekt C2) of the German Research Foundation
(DFG). The work of Viatcheslav Melas and Petr Shpilev was partly supported by Russian Foundation for Basic Research (project no. 17-01-00161).

\begin{small}
 \setlength{\bibsep}{4pt}

\end{small}

\end{document}